\documentclass[11pt]{article}
\usepackage{amssymb}
\usepackage{CJK}
\usepackage{amsmath,amsfonts,mathrsfs,amssymb}
\usepackage{indentfirst}

\numberwithin{equation}{section}

\begin{document}

\title{$\lambda$-Biharmonic  hypersurfaces in the product space $L^{m}\times \mathbb{R}$}
\date{}
\maketitle
\vspace{-1.5cm}
\centerline{\author{Chao Yang$^{\ast}$, Zhen Zhao}}
\begin{center}
\small{College of Mathematics and Statistics, Northwest Normal University, Lanzhou, 730070,
China}
\end{center}

\vskip.3cm

\begin{quotation}
{\bf Abstract:}
{In this paper, we study $\lambda$-biharmonic hypersurfaces in the product space $L^{m}\times\mathbb{R}$, where $L^{m}$ is an Einstein space and $\mathbb{R}$ is a real line. We prove that $\lambda$-biharmonic hypersurfaces with constant mean curvature in $L^{m}\times\mathbb{R}$ are either minimal or vertical cylinders, and obtain some classification results for $\lambda$-biharmonic hypersurfaces under various constraints. Furthermore, we investigate $\lambda$-biharmonic hypersurfaces in the product space $L^{m}(c)\times\mathbb{R}$, where $L^{m}(c)$ is a space form with constant sectional curvature $c$, and categorize hypersurfaces that are either totally umbilical or semi-parallel.}

\vskip.2cm
{\bf Keywords:}\ product space; $\lambda$-biharmonic  hypersurfaces; constant mean curvature; totally umbilical hypersurfaces; semi-parallel hypersurfaces

\vskip.2cm
{\bf Mathematics Subject Classification (2020)} 53C50

\end{quotation}

\footnote[0]{This work is supported by National Natural Science Foundation
of China (No. 12161078), Science and Technology Project of Gansu Province (No. 20JR5RA515).}

\footnote[0]{E-mail address: yangch@nwnu.edu.cn}

\section{Introduction}

Let $\varphi:M^{m}\rightarrow N^{m}\times\mathbb{R}$ be an isometric immersion. Denote by $E(\varphi)$ and $E_{2}(\varphi)$ the energy and bienergy functionals of $\varphi$, respectively. The functional
$$
E_{2,\lambda}(\varphi)=E_{2}(\varphi)+\lambda E(\varphi),
$$
is called the $\lambda$-bienergy functional of $\varphi$. We say that a map $\varphi$ is $\lambda$-biharmonic if it is a critical point of the functional $E_{2,\lambda}(\varphi)$. The Euler-lagrange equation of $E_{2,\lambda}(\varphi)$ yields
\begin{equation}\label{1.1}
\tau_{2}(\varphi)-\lambda \tau(\varphi)=0,
\end{equation}
where
\begin{align*}
&\quad \quad \quad \quad \ \ \tau(\varphi)=tr(\nabla d\varphi),\\
\tau_{2}(\varphi):=&\Delta\tau(\varphi)-tr_{g}\widetilde{R}(d\varphi,\tau(\varphi)d\varphi)\\
=&tr_{g}(\nabla^{\varphi}\nabla^{\varphi}\tau(\varphi)-\nabla^{\varphi}_{\nabla^{M}}\tau(\varphi))-tr_{g}\widetilde{R} (d\varphi,\tau(\varphi)d\varphi),
\end{align*}
$\nabla^{\varphi}$ and $\nabla$ are induced connections of vector bundle $\varphi^{-1}TN^{m+1}$ and Levi-Civita connection on $M^{m}$ respectively, and $\widetilde{R}$ is the curvature operator on $N^{m+1}$.

Special cases are discussed, such as when $\lambda=0$ leading to automatically biharmonic hypersurfaces. In the case of the ambient space $N^{m+1}$ having constant curvature $c$, the equation \eqref{1.1} simplifies to $\Delta \overrightarrow{H}=(mc-\lambda)\overrightarrow{H}$, where $\overrightarrow{H}$ is the mean curvature vector field, defining hypersurfaces with proper mean curvature vector fields.

The first result on $\lambda$-biharmonic hypersurfaces was obtained in \cite{Chen 1988-1} by B. Y. Chen in 1988. It was proved that every $\lambda$-biharmonic hypersurface in $\mathbb{R}^{3}$ is either minimal or an open portion of a circular cylinder. Furtherly, A. Ferr\'{a}ndez and P. Lucas showed in \cite{Ferrandez 1991} that $\lambda$-biharmonic hypersurfaces in $\mathbb{R}^{m+1}$ with at most two distinct principal curvatures are minimal or locally isometric to $\mathbb{R}^{k}\times S^{m-k}(a)$. For $\lambda$-biharmonic hypersurfaces with three or more distinct  principal curvatures, it is difficult to obtain classification results. But geometers conjectured that their mean curvatures are constants. In 1995, T. Hasanis and T. Vlachos showed that $\lambda$-biharmonic  hypersurfaces in $\mathbb{R}^{4}$ have constant mean curvatures in \cite{Hasanis 1995}.
In 2015, Y. Fu derived in \cite{Fu 2015-3} that $\lambda$-biharmonic hypersurface in $\mathbb{R}^{m+1}$ with at most three distinct principal curvatures have constant mean curvatures. Afterwards, Y. Fu and X. Zhan proved in \cite{Fu 2021} that $\lambda$-biharmonic hypersurfaces in $\mathbb{R}^{5}$ have constant mean curvatures.

For $\lambda$-biharmonic hypersurfaces in non-flat space forms, C. Yang, J.-C. Liu and L. Du gave in \cite{Yang 2024} classification results under the condition that hypersurfaces have at most two distinct principal curvatures. Scholars have also discussed $\lambda$-biharmonic  hypersurfaces in pseudo-Riemannian space forms (see \cite{Arvanitoyeorgos 2007, Arvanitoyeorgos 2009, Arvanitoyeorgos 2013, Du 2017, Du 2023, Ferrandez 1992-1, Liu 2014, Liu 2016, Liu 2017, Yang 2023, Yang 2024}), aiming to prove that their mean curvatures are constants.

Naturally, it is meaningfull to study $\lambda$-biharmonic hypersurfaces in the product space $L^{m}\times \mathbb{R}$ of an Einstein space $L^{m}$ and a real line $\mathbb{R}$.
The focus of this paper is on $\lambda$-biharmonic hypersurfaces in the product space $L^{m}\times\mathbb{R}$. Significant results are presented, showing that hypersurfaces with constant mean curvature in $L^{m}\times\mathbb{R}$ are either minimal or vertical cylinders. Moreover, classifications for hypersurfaces in $L^{m}(c)\times\mathbb{R}$ as totally umbilical or semi-parallel are discussed.

\section{Preliminaries}

Let $\varphi:M^{m}\rightarrow (N^{m}\times\mathbb{R},g^{N}+\text dt^{2})$ be a hypersurface with unit normal vector $\xi$. $T$ denotes the tangential component of $\partial_{t}$ along the tangent space to $M^{m}$ and $\alpha$ be the angle made by $\partial_{t}$ and $\xi$. We can decompose $\partial_{t}$ in the following form
$$
\partial_{t}=T+\cos\alpha\xi,
$$
where $\cos \alpha=\langle\partial_{t},\xi\rangle$. Let $\widetilde{Ric}$ be the Ricci curvature tensor of $N^{m}\times\mathbb{R}$,
$A$ and $H$ be the shape operator and the mean curvature of the hypersurface $M^{m}$, then the Laplacian of $\theta:=\cos \alpha$ is given by
\begin{equation}\label{2.1}
\Delta\theta=-m\langle\nabla H,\partial_{t}\rangle-\theta(|A|^{2}+\widetilde{Ric}(\xi,\xi)).
\end{equation}

Since $\partial_{t}$ is parallel in $N^{m}\times\mathbb{R}$, we have
\begin{equation}\label{2.2}
\left\{\begin{array}{l}
\nabla_{X}T=\cos\alpha AX\\
X(\cos\alpha)=-\langle AX,T\rangle,
\end{array}
\right.
\end{equation}
for every tangent vector field $X$ on $M^{m}$.

Let $\widetilde{Ric}$, $\widetilde{R}$ and $\widetilde{S}$ be the Ricci curvature tensor, Riemannian curvature tensor, and scalar curvature of $N^{m+1}\times\mathbb{R}$, respectively. Denote by  $Ric$, $R$ and $S$ corresponding ones of $M^{m}$. Then the relationships of Ricci curvature tensors and scalar curvatures between the hypersurface and the ambient space are given by
$$
\widetilde{Ric}(X,Y)=Ric(X,Y)+\langle AX,AY\rangle-mH\langle AX,Y\rangle+\widetilde{R}(X,\xi,Y,\xi),
$$
and
\begin{equation}\label{2.3}
\widetilde{S}=S+|A|^{2}-m^{2}H^{2}+2\widetilde{Ric}(\xi,\xi).
\end{equation}
The second fundamental form $B$ and the shape operator $A$ satisfy
$$
\langle B(X,Y),\xi\rangle=\langle A(X),Y\rangle.
$$
Moreover, the Codazzi equation is given by
$$
(\widetilde{R}(X,Y)Z)^{\bot}=(\nabla_{X}B)(Y,Z)-(\nabla_{Y}B)(X,Z).
$$

If the hypersurface $M^{m}$ is totally umbilical if and only if
$$
\langle B(X,Y),\xi\rangle=\langle X,Y\rangle H.
$$
We say the hypersurface $M^{m}$ is semi-parallel if and only if
$$
B(R(X,Y)U,V)+B(U,R(X,Y)V)=0,
$$
for any $X,Y,U,V\in \Gamma (M)$.

In this paper, we denote by $L^{m}$ and $L^{m}(c)$ the  $m$-dimensional Einstein manifold and the space with constant curvature $c=1,-1$ or $0$, respectively. For simplicity, we call $L^{m}\times\mathbb{R}$ Einstein product space, and $L^{m}(c)\times\mathbb{R}$ constant curvature product space. 
The Riemannian curvature tensor of $L^{m}(c)\times\mathbb{R}$ is given by (\cite{Hasanis 1995})
\begin{equation*}
\begin{aligned}
\widetilde{R}(X,Y)Z=&c\{\langle Y,Z\rangle X-\langle X,Z\rangle Y-\langle Y,\partial_{t}\rangle\langle Z,\partial_{t}\rangle X+
\langle X,\partial_{t}\rangle\langle Z,\partial_{t}\rangle Y\\&+\langle X,Z\rangle\langle Y,\partial_{t}\rangle\partial_{t}-
\langle Y,Z\rangle\langle X,\partial_{t}\rangle\partial_{t}\},
\end{aligned}
\end{equation*}
where $X,Y,Z$ are vector fields on $L^{m}(c)\times\mathbb{R}$.

Recall in \cite{Dillen 2009}, that a rotation hypersurface $M^{m}$ in the product space $S^{m}\times\mathbb{R}$, parametrizing the profile curve as
$$
\gamma(s)=(\cos s,0,\cdots,0,\sin s,h(s)),
$$
for some smooth function $h$, then the rotation hypersurface $M^{m}$ can be parametrized as
$$
f(s,v_{1},\cdots,v_{m-1})=(\cos s,\varphi_{1}(v_{1},\cdots,v_{m-1})\sin s,\cdots,\varphi_{m}(v_{1},\cdots,v_{m-1})\sin s,h(s)),
$$
where $\varphi=(\varphi_{1},\cdots\varphi_{m})$ is an orthogonal parametrization of the unit sphere $S^{m-1}$ in $\mathbb{R}^{m}$,
we compute the principal curvatures as follows
$$
\lambda_{1}=-\frac{h^{\prime\prime}(s)}{(1+h^{\prime}(s)^{2})^{\frac{3}{2}}},\ \
\lambda_{2}=-\frac{h^{\prime}(s)\cot s}{(1+h^{\prime}(s)^{2})^{\frac{1}{2}}}.
$$
The similar argument for a rotation hypersurface $M^{m}$ in a  product space $H^{m}\times\mathbb{R}$,
the principal curvatures as follows
$$
\lambda_{1}=-\frac{h^{\prime\prime}(s)}{(1+h^{\prime}(s)^{2})^{\frac{3}{2}}},\ \
\lambda_{2}=-\frac{h^{\prime}(s)\tan a(s)}{(1+h^{\prime}(s)^{2})^{\frac{1}{2}}}.
$$

To obtain main results, the following lemmas are required.

\section{Some Lemmas}

\vskip.2cm
{\bf Lemma 3.1}\quad \emph{Let $\varphi:(M^{m},g)\rightarrow (N^{m}\times\mathbb{R}, g^{N}+\text{d}t^{2})$ be an isometric immersion.  Then $\varphi$ is $\lambda$-biharmonic if and only if both its component maps $\pi_{1}\circ\varphi:(M,g)\rightarrow (N,g^{N})$ and $\pi_{2}\circ\varphi:(M,g)\rightarrow(\mathbb{R},\text{d}t^{2})$ are $\lambda$-biharmonic maps. In particular, the height function $h=\pi_{2}\circ\varphi$ of a $\lambda$-biharmonic hypersurface is a $\lambda$-biharmonic function on the hypersurface.}

\vskip.2cm
{\bf Proof}\quad Since
$$
\text{d}\varphi(X)=\text{d}(\pi_{1}\circ\varphi)(X)+\text{d}(\pi_{2}\circ\varphi)(X), \forall X\in\Gamma M,
$$
we know
$$
\nabla^{\varphi}_{X}(\text{d}\varphi(Y))=\nabla^{\varphi}_{X}(\text{d}(\pi_{1}\circ\varphi))(Y)+\nabla^{\varphi}_{X}(\text{ d}(\pi_{2}\circ\varphi))(Y),\forall X,Y\in\Gamma M.
$$
It follows that
\begin{align*}
\tau(\varphi)=&\sum^{m}_{i=1}\{\nabla^{\varphi}_{e_{i}}(\text{d}\varphi(e_{i}))-\text{ d}\varphi(\nabla_{e_{i}}{e_{i}})\}\\
             =&\sum^{m}_{i=1}\{\nabla_{e_{i}}^{\pi_{1}\circ\varphi}(\text{d}(\pi_{1}\circ\varphi)(e_{i}))-\text{d}(\pi_{1}\circ\varphi)(\nabla^{M}_{e_{i}}e_{i})\}\\&+
\sum^{m}_{i=1}\{\nabla_{e_{i}}^{\pi_{2}\circ\varphi}\text{d}(\pi_{2}\circ\varphi)(e_{i})-\text{ d}(\pi_{2}\circ\varphi)(\nabla^{M}_{e_{i}}e_{i})\}\\
             =&\tau(\pi_{1}\circ\varphi)+\tau(\pi_{2}\circ\varphi).
\end{align*}
Since
$$
\tau_{2}(\varphi)=-J^{\varphi}(\tau(\varphi)),
$$
$$
J^{\varphi}(X):=-\{tr(\nabla^{\varphi}\nabla^{\varphi}-\nabla^{\varphi}_{\nabla^{M}})X-tr\widetilde{R}(\text{d}\varphi,X)\text{d}\varphi\},
$$
we get
\begin{align*}
\tau_{2}(\varphi)&=-J^{\varphi}(\tau(\varphi))\\
&=- J^{\pi_{1}\circ\varphi}(\tau(\pi_{1}\circ\varphi))- J^{\pi_{2}\circ\varphi}(\tau(\pi_{2}\circ\varphi))\\
&=\tau_{2}(\pi_{1}\circ\varphi)+\tau_{2}(\pi_{2}\circ\varphi).
\end{align*}
That $\varphi$ is $\lambda$-biharmonic means that $\tau_{2}(\varphi)-\lambda \tau(\varphi)=0$, i.e.
$$
\left\{\begin{array}{l}
\tau_{2}(\pi_{1}\circ\varphi)-\lambda\tau(\pi_{1}\circ\varphi)=0,\\
\tau_{2}(\pi_{2}\circ\varphi)-\lambda\tau(\pi_{2}\circ\varphi)=0.
\end{array}
\right.
$$
Then we obtain the lemma.
\hfill$\square$

\vskip.2cm
{\bf Lemma 3.2}\quad \emph{Let $\varphi:(M^{m},g)\rightarrow (N^{m}\times\mathbb{R},g^{N}+\text dt^{2})$ be an isometric immersion. If the height function $h=\pi_{2}\circ\varphi$ is $\lambda$-biharmonic, then we have
$$
\Delta^{2}h=\lambda\Delta h.
$$
}

\vskip.2cm
{\bf Proof}\quad Calculate $\tau(h)$, we derive
\begin{align*}
\tau(h)&=tr(\nabla \text{d}h)\\
&=\nabla_{e_{i}}(\text{d}h(e_{i}))-\text{d}h(\nabla_{e_{i}}e_{i})\\
&=e_{i}e_{i}h-(\nabla_{e_{i}}e_{i})h\\
&=\Delta h.
\end{align*}
Considering the curvature operator $\bar{R}$ on Euclidean space is zero, we know that
$$
\tau_{2}(h)=\Delta\tau(h)-tr\bar{R}(\text{d}h,\tau(h)dh)=\Delta^{2}h.
$$
So, $\tau_{2}(h)=\lambda \tau(h)$ is equivalent to $\Delta^{2}h=\lambda\Delta h$.
\hfill$\square$

\vskip.2cm
{\bf Lemma 3.3}\quad \emph{Let $M^{m}$ is a $\lambda$-biharmonic hypersurface in $N^{m}\times\mathbb{R}$. Then we have the following identity
\begin{equation}\label{3.1}
\Delta(H\theta)=\lambda H\theta.
\end{equation}
}

\vskip.2cm
{\bf Proof}\quad Let $h$ be the height function of the hypersurface, we know from \cite{3} that  $\Delta h=m\theta H$. Furthermore
$$
\Delta^{2}h=\Delta(\Delta h)=m\Delta(H\theta),
$$
which together with $\Delta^{2}h=\lambda\Delta h$ (cf. Lemma 3.2) deduced the lemma.
\hfill$\square$

\vskip.2cm
{\bf Lemma 3.4(\cite{Yau 1976})}\quad \emph{(Yau's maximum principle)
\begin{itemize}
  \item Let $M^{m}$ is a complete Riemannian manifold with nonnegative Ricci curvature and $u$ is a positive smooth harmonic function. Then $u$ is a constant function;
  \item Let $M^{m}$ is a complete Riemannian manifold and $u$ is a nonnegative smooth subharmonic function. Then $\int_{M}u^{p}=+\infty$ for $p>1$, unless $u$ is a constant function.
\end{itemize}}

\vskip.2cm
{\bf Lemma 3.5(\cite{Luo 2017})}\quad \emph{Let $u\in(0,C](C>0)$ be a superharmonic function on complete noncompact manifold $M^{m}$. If
$$
\int_{M}(\log^{(k)}\frac{Ce^{(k)}}{u})^{p}\text{d}v_{g}<+\infty, \text for  \text some p>0, k\in\mathbb{N},
$$
where $\log^{(1)}=\log,e^{(1)}=e,\log^{(k)}=\log(\log^{(k-1)}),e^{(k)}=e^{e^{(k-1)}}$, then $u$ is a constant.}

\vskip.2cm
{\bf Lemma 3.6}\quad \emph{A hypersurface $\varphi :M^{m}\rightarrow N^{m+1}$ in a Riemannian manifold is $\lambda$-biharmonic if and only if
\begin{equation}\label{3.2}
\left\{\begin{array}{l}
\Delta H-H|A|^{2}+H\widetilde{Ric}(\xi,\xi)-\lambda H=0,\\
2A(\nabla H)+\frac{m}{2}\nabla H^{2}-2H(\widetilde{Ric}(\xi))^{T}=0,
\end{array}
\right.
\end{equation}
where $\widetilde{Ric}$ denotes Ricci curvature operator on $N^{m+1}$, and $\langle \widetilde{Ric}(Z),W\rangle:=\widetilde{Ric}(Z,W)$ for any $Z,W\in\Gamma M$.
}

\vskip.2cm
{\bf Proof}\quad Choose a local orthonormal frame $\{e_{i}\}, i=1,\cdots,m$ on $M^{m}$, such that  $\{d\varphi(e_{1}),\cdots,d\varphi(e_{m}),\xi\}$ is an local orthonormal frame of $N^{m+1}$. Notice that $\tau(\varphi)=mH\xi$, we compute the bitension field of $\varphi$
\begin{align*}
\tau_{2}(\varphi)&=\sum_{i=1}^{m}\{\nabla^{\varphi}_{e_{i}}\nabla^{\varphi}_{e_{i}}(mH\xi)-\nabla^{\varphi}_{\nabla_{e_{i}}{e_{i}}}(mH\xi)-\widetilde{R}(\text{d}\varphi(e_{i}),mH\xi)\text{d}\varphi(e_{i})\}\\
&=m(\Delta H)\xi-2mA(\nabla H)-mH\Delta^{\varphi}\xi-mH\sum_{i=1}^{m}\widetilde{R}(\text{d}\varphi(e_{i}),\xi)\text{d}\varphi(e_{i}),
\end{align*}
where $\tilde{R}$ denotes the Riemannian curvature operator of the ambient space $N^{m+1}$.

As
$$
\sum^{m}_{i=1}\langle\widetilde{R}(\text{d}\varphi(e_{i}),\tau(\varphi))\text{d}\varphi(e_{i}),\xi\rangle=-mH\widetilde{Ric}(\xi,\xi),
$$
and
\begin{align*}
\langle\Delta^{\varphi}\xi,\xi\rangle&=\sum^{m}_{i=1}\langle\widetilde{\nabla}_{e_{i}}\widetilde{\nabla}_{e_{i}}\xi-\widetilde{\nabla}_{\nabla_{e_{i}}{e_{i}}}\xi,\xi\rangle\\
&=\sum^{m}_{i=1}\langle\widetilde{\nabla}_{e_{i}}\xi,\widetilde{\nabla}_{e_{i}}\xi\rangle=|A|^{2},
\end{align*}
the normal part of $\tau_{2}(\varphi)-\lambda\tau(\varphi)$ can be expressed as
\begin{equation*}
\begin{aligned}
(\tau_{2}(\varphi)-\lambda\tau(\varphi))^{\bot}&=\langle\tau_{2}(\varphi)-\lambda\tau(\varphi),\xi\rangle\xi\\
&=\Delta H-H|A|^{2}+H\widetilde{Ric}(\xi,\xi)-\lambda H.
\end{aligned}
\end{equation*}

Because of
\begin{align*}
\langle\Delta^{\varphi}\xi,e_{k}\rangle e_{k} &=\sum^{m}_{i,k=1}\langle\widetilde{\nabla}_{e_{i}}\widetilde{\nabla}_{e_{i}}\xi-\widetilde{\nabla}_{\nabla_{e_{i}}{e_{i}}}\xi,e_{k}\rangle e_{k}\\
&=m(\nabla H)-(\widetilde{Ric}(\xi))^{\top},
\end{align*}
and
$$
\sum^{m}_{i,k=1}\langle\widetilde{R}(d\varphi(e_{i}),\xi)d\varphi(e_{i}),e_{k}\rangle e_{k}=-[\widetilde{Ric}(\xi,e_{k})]e_{k}=-(\widetilde{Ric}(\xi))^{\top},
$$
the tangent part of $\tau_{2}(\varphi)-\lambda\tau(\varphi)$ is
\begin{align*}
(\tau_{2}(\varphi)-\lambda\tau(\varphi))^\top&=\langle\tau_{2}(\varphi)-\lambda\tau(\varphi),e_{k}\rangle e_{k}\\
&=-2A(\nabla H)-\frac{m}{2}(\nabla H^{2})+2H(\widetilde{Ric}(\xi))^{\top}).
\end{align*}
The hypersurface $M^{m}$ is $\lambda$-biharmonic if and only if the tangent and normal parts of $\tau_{2}(\varphi)-\lambda\tau(\varphi)$ are equal to zero, i.e. \eqref{3.2} formula holds.
\hfill$\square$

\section{$\lambda$-biharmonic hypersurface in $L^{m}\times\mathbb{R}$}

Using Lemmas in Section 3, we can give some classification results for $\lambda$-biharmonic hypersurface in a product space $L^{m}\times\mathbb{R}$ of an Einstein space $L^{m}$ and a real line $\mathbb{R}$. In the following, we classify the $\lambda$-biharmonic hypersurface with constant mean curvature in such a product space.

\vskip.2cm
{\bf Theorem 4.1}\quad \emph{Let $L^{m}$ be an Einstein manifold and $M^{m}$ is a $\lambda$-biharmonic $(\lambda\geq0)$ hypersurface in $L^{m}\times\mathbb{R}$ with constant mean curvature. Then $M^{m}$ is minimal, or a vertical cylinder over a $\lambda$-biharmonic hypersurface in $L^{m}$.}

\vskip.2cm
{\bf Proof}\quad If $H=0$, then $M^{m}$ is minimal.

If $H\neq0$, then \eqref{3.1} can be reduced to
\begin{equation}\label{4.1}
\Delta\theta=\lambda \theta.
\end{equation}
It is easy to find that
\begin{equation}\label{4.2}
\widetilde{Ric}(\xi,\xi)=\mu(1-\theta^{2}),
\end{equation}
where $\mu$ is a constant.

Substitute \eqref{4.2} into \eqref{2.1} and the first equation of \eqref{3.2}, and considering $H$ is a constant, we have
\begin{equation}\label{4.3}
\Delta\theta=-\theta(2|A|^{2}+\lambda).
\end{equation}
Combine \eqref{4.1} and \eqref{4.3}, we have $\theta(|A|^{2}+\lambda)=0$.
As $|A|^{2}>0$, $\lambda\geq0$, the equation $\theta(|A|^{2}+\lambda)=0$ implies $\theta\equiv0$, i.e. $\partial_{t}$ is tangent to the hypersurface. So, we conclude that the hypersurface $M^{m}$ is a vertical cylinder over a $\lambda$-biharmonic hypersurface in $L^{m}$.
\hfill$\square$

\vskip.2cm
{\bf Theorem 4.2}\quad \emph{Let $L^{m}$ be an Einstein manifold and $M^{m}$ is a complete $\lambda$-biharmonic $(\lambda\geq0)$ hypersurface in $L^{m}\times\mathbb{R}$ with constant angle function. If the mean curvature $H$ is nonnegative and $H\in L^{p}(M),1<p<\infty$, then $M^{m}$ is minimal, or a vertical cylinder over a $\lambda$-biharmonic hypersurface in $L^{m}$.}

\vskip.2cm
{\bf Proof}\quad If $\theta=\langle\xi,\partial_{t}\rangle\equiv0$, then the hypersurface $M^{m}$ is a vertical cylinder over a $\lambda$-biharmonic hypersurface in $L^{m}$.

If the constant $\theta\neq0$, then, we have $\Delta H=\lambda H$ by \eqref{3.1}. As $\lambda\geq0$ and $H$ is nonnegative, it follows that $\Delta H\geq0$. Applying $\text Yau^{,}s$ maximum principle, we find that $H$ is a constant. We complete the proof by using Theorem 4.1.
\hfill$\square$

\vskip.2cm
{\bf Theorem 4.3}\quad \emph{Let $L^{m}$ be an Einstein manifold and $M^{m}$ is a totally umbilical $\lambda$-biharmonic $(\lambda\geq0)$ hypersurface in $L^{m}\times\mathbb{R}$ with constant angle function. Then $M^{m}$ is minimal, or a vertical cylinder over a $\lambda$-biharmonic hypersurface in $L^{m}$.}

\vskip.2cm
{\bf Proof}\quad Choose $\{e_{1},e_{2},\cdots e_{m}\}$ be a local orthonormal frame on $M^{m}$, then  $A(e_{i})=He_{i},i=1,2,\cdots m$ and $|A|^{2}=mH^{2}$.

Substitute $\widetilde{Ric}(\xi,\xi)=\mu(1-\theta^{2})$ into the first equation of \eqref{3.2}, we have
$$
\Delta H-mH^{3}+H\mu(1-\theta^{2})-\lambda H=0.
$$
If $\theta\neq0$, we know $\Delta H=\lambda H$ from \eqref{3.1}. So the above equation can be reduced to
$$
mH^{3}=H\mu(1-\theta^{2}),
$$
which means that $H$ is a constant. Applying Theorem 4.1, we finish the proof.
\hfill$\square$

\vskip.2cm
{\bf Proposition 4.4}\quad \emph{Let $M^{m}$ is a complete $\lambda$-biharmonic $(\lambda\geq0)$ hypersurface with non-negative Ricci curvature in $L^{m}\times\mathbb{R}$ . Assume that
$$
\int_{M}H^{2p}dv_{g}<+\infty,
$$
and
$$
\int_{M}(\log^{(k)}\frac{e^{(k)}}{\theta^{2}+\varepsilon})^{q}dv_{g}<+\infty,
$$
for some $p>1, q>0,k\in\mathbb{N}$ and $\varepsilon>0$. Then $M^{m}$ is minimal, or a vertical cylinder over a $\lambda$-biharmonic hypersurface in $L^{m}$.}

\vskip.2cm
{\bf Proof}\quad It follows from Lemma 3.3 that
$$
\Delta(H\theta)^{2}=2|\nabla(H\theta)|^{2}+2\lambda(H\theta)^{2}\geq0.
$$
As $-1\leq\theta\leq1$, we have
$$
\int_{M}{(H\theta)^{2}}^{p}dv_{g}\leq\int_{M}H^{2p}dv_{g}<+\infty, p>1,
$$
Then by Yau's maximum principle, we conclude that $H\theta$ is a constant.

As $\Delta h=mH\theta$ and $H\theta$ is a constant, the Ricci identity
$$
\Delta\nabla_{i}h=\nabla_{i}\Delta h+Ric_{ij}\nabla_{j}h.
$$
can be simplify as
$$
\Delta\nabla_{i}h=Ric_{ij}\nabla_{j}h.
$$
From $h=\pi_{2}\circ\varphi$, we have
$$
\langle\Delta T,T\rangle=Ric(T,T).
$$

Due to
\begin{equation}\label{4.4}
\frac{1}{2}\Delta|T|^{2}=|\nabla T|^{2}+\langle\Delta T,T\rangle=|\nabla T|^{2}+Ric(T,T)\geq0.
\end{equation}
Also because of $\partial_{t}=T+\theta\xi$, we have $\langle T,T\rangle=1-\theta^{2}$.
So
\begin{equation}\label{4.5}
\frac{1}{2}\Delta|T|^{2}=\frac{1}{2}\Delta(1-\theta^{2})=-\frac{1}{2}\Delta(1+\theta^{2}).
\end{equation}
Combining \eqref{4.4} and \eqref{4.5}, we have
$$
-\frac{1}{2}\Delta (\theta^{2}+1)\geq0.
$$
Therefore $\Delta(\theta^{2}+1)\leq0$, that is $\theta^{2}+\varepsilon$ is a superharmonic function. From Lemma 3.5, we obtain that $\theta^{2}+1$ is a constant, i.e. $\theta$ is a constant. It together with that $H\theta$ is a constant imply that $H$ is also a constant. By Theorem 4.1, the proof is completed.
\hfill$\square$

\vskip.2cm
{\bf Proposition 4.5}\quad \emph{let $M^{m}$ be a complete $\lambda$-biharmonic hypersurface with non-negative Ricci curvature in a product space $L^{m}\times \mathbb{R}$. Assume that
\begin{itemize}
  \item $H$ is harmonic and bounded from below, or
  \item $\theta$ is harmonic and scalar curvature of $M^{m}$ is a constant.
\end{itemize}
Then $M^{m}$ is either minimal, or a vertical cylinder over a $\lambda$-biharmonic hypersurface in $L^{m}$.}

\vskip.2cm
{\bf Proof}\quad (i) Assume that $\Delta H=0$ and $H\geq-C$ for some positive $C$, take $u=H+C+\varepsilon>0$, where $\varepsilon$ is some positive constant. Then $\Delta u=\Delta H=0$. According to Yau's maximum principle, we know $u$ is a constant. Hence, $H=u-C-\varepsilon$ is a constant. By Theorem 4.1, the result follows.

(ii)We suppose that $\Delta\theta=0$ and the scalar curvature $S$ is a constant. Set $u=\theta+2$, then $u>0$ and $\Delta u=\Delta \theta=0$. It follows from Yau's maximum principle that $u$ and $\theta$ are constant.

If $\theta\equiv0$, then the $\lambda$-biharmonic hypersurface is a vertical cylinder over a $\lambda$-biharmonic hypersurface in $L^{m}$.

If $\theta\neq0$, we have $\Delta H=\lambda H$. Then using the equation $\widetilde{Ric}(\xi,\xi)=\mu(1-\theta^{2})$, the first equation of \eqref{3.2} use can be rewritten as
$$
H[|A|^{2}-\mu(1-\theta^{2})]=0.
$$
When $H\equiv0$, then $M^{m}$ is minimal.
When $H\neq0$ at some point $p\in M$, then the above formula tells us that
$$
|A|^{2}=\mu(1-\theta^{2}).
$$
By \eqref{2.3}, one get $S=\widetilde{S}-3\mu(1-\theta^{2})+m^{2}H^{2}$. Noting that $\widetilde{S}=\mu m$. It follows that $\nabla S=2m^{2}H\nabla H$. We find from $S$ is a constant $\nabla S=0$, which together with $2m^{2}H\nabla H=\nabla S$ deduced that $H$ is a constant.

Next we can use Theorem 4.1 to complete the proof.

In particular, when the product space is $L^{m}(c)\times \mathbb{R}$, we have the following Theorem.
\hfill$\square$

\section{$\lambda$-biharmonic hypersurface in $L^{m}(c)\times\mathbb{R}$}

\vskip.2cm
{\bf Lemma 5.1}\quad \emph{A hypersurface $M^{m}$ in a product space $L^{m}(c)\times\mathbb{R}$ is $\lambda$-biharmonic if and only if
\begin{equation}\label{5.1}
\left\{\begin{array}{l}
\Delta H-H[|A|^{2}-c(m-1)\sin^{2}\alpha+\lambda]=0,\\
A(\nabla H)+\frac{m}{2} H\nabla H+c(m-1)\cos\alpha HT=0.
\end{array}
\right.
\end{equation}}

\vskip.2cm
{\bf Proof}\quad Set $\{e_{i}\}, i=1,\cdots,m$ be a local orthonormal frame on $M^{m}$. Using \eqref{2.3}, we get

$$
\widetilde{Ric}(\xi,\xi)=\sum^{m}_{i=1}\langle\widetilde{R}(e_{i},\xi)\xi,e_{i}\rangle=c(m-1)\sin^{2}\alpha
$$
and
$$
(\widetilde{Ric}(\xi))^{T}=\sum^{m}_{i=1}\langle\widetilde{R}(e_{i},\xi)e_{k},e_{i}\rangle e_{k}=-c(m-1)\cos\alpha T.
$$
Substituting the above two equations into $\lambda$-biharmonic equations \eqref{3.2}, we obtain the lemma.
\hfill$\square$

\vskip.2cm
{\bf Theorem 5.2}\quad \emph{A rotation hypersurface $M^{m}$ in $L^{m}(c)\times\mathbb{R}$ is $\lambda$-biharmonic, then we have
\begin{equation}\label{5.2}
(\frac{m}{2}H-\alpha^{\prime}(s)\cos \alpha)H^{\prime}-c(m-1)\sin \alpha H=0.
\end{equation}
}

\vskip.2cm
{\bf Proof}\quad As in \cite{Dillen 2009}, we choose a local orthonormal frame
 $\{e_{1},e_{2},\cdots,e_{m}\}$, so that
$$
\nabla H=e_{1}(H)e_{1}=-\cos \alpha H^{\prime},
$$
$$
A(e_{1})=\lambda_{1} e_{1}, A(e_{i})=\lambda_{i}e_{i},   2\leq i\leq m.
$$
Combining $\langle T,T\rangle=\sin^{2}\alpha$, the second equation of \eqref{5.1} can be rewriten as
$$
(\frac{m}{2}H+\lambda_{1})H^{\prime}-c(m-1)\sin\alpha H=0.
$$
\hfill$\square$

\vskip.2cm
{\bf Theorem 5.3}\quad \emph{Any totally umbilical $\lambda$-biharmonic ($\lambda\geq0$) hypersurface $M^{m}$ in $L^{m}(c)\times \mathbb{R}$ is minimal.}

\vskip.2cm
{\bf Proof}\quad Since $M^{m}$ is a totally umbilical $\lambda$-biharmonic hypersurface, we have
 $|A|^{2}=mH^{2}$ and $A(\nabla H)=H\nabla H$. Then the equation \eqref{5.1} can be rewritten as
\begin{equation}\label{5.3}
\left\{\begin{array}{l}
\Delta H-H[mH^{2}-c(m-1)\sin^{2}\alpha+\lambda]=0\\
\frac{m+2}{2}\nabla H+c(m-1)\cos\alpha T=0.
\end{array}
\right.
\end{equation}

If $H\equiv0$, then $M^{m}$ is minimal.

Now, we assume that $H\neq0$ on an open set $\Omega$.

(i) When $\sin\alpha\equiv0$, we know $\cos\alpha=\pm1$. Take $X=\nabla H$ in the second equation of \eqref{2.2}, and combine $A(\nabla H)=H\nabla H$ and the second equation of \eqref{5.3}, we have $H|\nabla H|^{2}=0$, which implies that $|\nabla H|=0$, i.e. $H$ is non-zero constant. Thus, we derive from the first equation of \eqref{5.3} that $\lambda=-mH^{2}$, which contradicts with $\lambda\geq0$.

(ii) When $\sin\alpha\neq0$ at some point $p\in M^{m}$. Considering $|T|^{2}=\sin^{2}\alpha$, we can choose a local orthonormal frame  $\{e_{1},e_{2},\cdots,e_{m}\}$ on $M^{m}$, such that
\begin{equation}\label{5.4}
T=\sin\alpha e_{1}.
\end{equation}
Combining $\nabla H=\sum_{i}e_{i}(H)e_{i}$ and \eqref{5.4}, we obtain from the second equation of \eqref{5.3} that
$$
e_{2}(H)=e_{3}(H)=\cdots=e_{m}(H)=0,
$$
and
\begin{equation}\label{5.5}
e_{1}(H)=-\frac{2c(m-1)}{m+2}\sin\alpha\cos\alpha.
\end{equation}
Using the first equation of \eqref{2.2}, we have
$$
\langle \cos Ae_{j},e_{j}\rangle=H\cos\alpha=\langle\nabla_{e_{j}}T,e_{j}\rangle
=e_{j}\langle T,e_{j}\rangle-\langle T,\nabla_{e_{j}}e_{j}\rangle=-\sin\alpha\langle e_{1},\nabla_{e_{j}}e_{j}\rangle,
$$
which deduce that
$$
\langle e_{1},\nabla_{e_{j}}e_{j}\rangle=-H\cot\alpha.
$$
Then we have
\begin{equation}\label{5.6}
\begin{aligned}
\Delta H&=e_{1}e_{1}(H)+\sum_{i}(\nabla_{e_{i}}e_{i})(H)\\
&=e_{1}e_{1}(H)+(m-1)\cot\alpha He_{1}(H).
\end{aligned}
\end{equation}
Put $X=e_{1}$ in the second equation of \eqref{2.2}, and combine $T=\sin\alpha e_{1}$, we obtain
\begin{equation}\label{5.7}
H=e_{1}(\alpha).
\end{equation}
Differentiating \eqref{5.5} along $e_{1}$, we derive
\begin{equation}\label{5.8}
e_{1}e_{1}(H)=-\frac{2c(m-1)}{m+2}\cos(2\alpha)H.
\end{equation}
Substitute \eqref{5.6} into the first equation of \eqref{5.3}, and combine \eqref{5.7} and \eqref{5.8}, we have
\begin{equation}\label{5.9}
-e_{1}(\alpha)(\frac{2c(m-1)}{m+2}\cos(2\alpha)+\frac{2c(m-1)^{2}}{m+2}\cos^{2}\alpha-c(m-1)\sin^{2}\alpha+m(e_{1}(\alpha))^{2}+\lambda)=0.
\end{equation}
If $e_{1}(\alpha)\equiv0$, then $H=e_{1}(\alpha)\equiv0$, a contradiction.

If $e_{1}(\alpha)\neq0$ on some neighborhood $\Omega$ on $M^{m}$, then the equation \eqref{5.9} reads
\begin{equation}\label{5.10}
\frac{2c(m-1)}{m+2}\cos(2\alpha)+\frac{2c(m-1)^{2}}{m+2}\cos^{2}\alpha-c(m-1)\sin^{2}\alpha+m(e_{1}(\alpha))^{2}+\lambda=0.
\end{equation}
It is shown (\cite{Calvaruso 2010}) that the Sine-Gordon equation
\begin{equation}\label{5.11}
e_{1}e_{1}(2\alpha)+c\sin(2\alpha)=0.
\end{equation}
Differentiating \eqref{5.10} along $e_{1}$ and combining \eqref{5.11}, we have
\begin{equation}\label{5.12}
2c\sin\alpha\cos\alpha[\frac{4(m-1)+2(m-1)^{2}+(m-1)(m+2)+m(m+2)}{m+2}]=0,
\end{equation}
If $\sin\alpha=1$, then $\cos\alpha=0$, we have $\alpha$ is a constant, so $e_{1}(\alpha)=0$, a contradiction.
Hence, \eqref{5.12} is equivalent to
$$
4m^{2}+3m-4=0,
$$
which implies that $m$ is a non-positive integer, a contradiction.
\hfill$\square$

\vskip.2cm
{\bf Theorem 5.4}\quad \emph{Any semi-parallel $\lambda$-biharmonic ($\lambda\geq0$) hypersurface in $S^{m}\times\mathbb{R} (m\geq3)$ is minimal or a vertical cylinder over a $\lambda$-biharmonic hypersurface in $S^{m}$.}

\vskip.2cm
{\bf Proof}\quad Let $M^{m}$ is a semi-parallel hypersurface in $S^{m}\times\mathbb{R} (m\geq3)$, then according to \cite{Van der Veken 2008}, one of the following is true,

(I) $M^{m}$ is totally umbilical;

(II) $M^{m}$ is an open part of rotation hypersurface with $\lambda_{1}\lambda_{2}=-\cos^{2}\alpha$, which $\lambda_{1}$ and $\lambda_{2}$ \\are principal curvatures in hypersurface $M^{m}$;

(III) $M^{m}\subset \tilde{M}^{m-1}\times\mathbb{R}$, where $\tilde{M}$ is a semi-parallel hypersurface of $S^{m}$.

Suppose $M^{m}$ satisfies  $\lambda_{1}\lambda_{2}=-\cos^{2}\alpha$, then we have
$$
\lambda_{1}=-\alpha^{\prime}(s)\cos\alpha,\ \
\lambda_{2}=-\sin\alpha\cot s.
$$
If $u=-\sin\alpha$, then
$$
uu^{\prime}\cot s=u^{2}-1.
$$
Solving this equation, we obtain $u=\pm\sqrt{1+C\sec^{2}(s)}$. However, it does not satisfy the equation \eqref{5.2}, which is a contradiction.
\hfill$\square$

Similarly, we can give the following result.

\vskip.2cm
{\bf Theorem 5.5}\quad \emph{Any semi-parallel $\lambda$-biharmonic ($\lambda\geq0$) hypersurface in $H^{m}\times\mathbb{R} (m\geq3)$ is minimal or a vertical cylinder over a $\lambda$-biharmonic hypersurface in $H^{m}$.}

\end{document}